\documentclass[10pt]{article}

\usepackage{graphicx}
\usepackage{indentfirst,csquotes}
\usepackage{amssymb,amsthm,amsmath}
\usepackage{xcolor,paralist,hyperref,titlesec,fancyhdr,etoolbox}

\titleformat{\section}
  {\normalfont\large\bfseries\centering}
  {\thesection}{1em}{}
\titlespacing*{\section}{0pt}{2ex}{1ex}

\hypersetup{
  colorlinks=true,
  linkcolor=black,
  filecolor=black,
  urlcolor=black
}

\title{4-D Visualization of Minkowski Quaternionic Point Set Operations}

\author{
Jakub \v{R}ada\\
\small Czech Technical University in Prague, Faculty of Architecture\\
\small Th\'akurova 9, 166~34 Prague~6 -- Dejvice, Czech Republic\\
\small \texttt{jakub.rada@cvut.cz}
\and
Daniela Velichov\'a\\
\small Slovak University of Technology in Bratislava, Faculty of Mechanical Engineering\\
\small N\'am.\ slobody~17, 812~31 Bratislava, Slovakia\\
\small \texttt{daniela.velichova@stuba.sk}
\and
Michal Zamboj\thanks{Corresponding author: michal.zamboj@pedf.cuni.cz}\\
\small Charles University, Faculty of Education\\
\small M.\ Rettigov\'e~4, 116~39 Prague~1, Czech Republic\\
\small \texttt{michal.zamboj@pedf.cuni.cz}
}

\date{\today}

\begin{document}

\maketitle

\begin{abstract}
The contribution emphasizes the geometric modeling point of view on Minkowski point set operations. In this paper, the Minkowski product is specified as the quaternionic product. Selected point sets are visualized using double orthogonal projection and perspective projection from four-dimensional to three-dimensional space. In particular, we demonstrate the generation of sets containing circles (Clifford torus, 3-sphere), lines (quadratic cone), or both.
\end{abstract}



\section{Introduction}
Minkowski point set operations provide a natural way to make geometric objects obey algebraic rules.  Two fundamental operations are the Minkowski sum $\oplus$ and the Minkowski product $\otimes$, defined for points represented by their position vectors. Let $A$ and $B$ be two point sets in $\mathbb{R}^n$. Then 
\begin{equation}
\label{eq:sum}
    A \oplus B := \{\, a + b \mid a \in A,\; b \in B \,\},
\end{equation}
\begin{equation}
\label{eq:product}
    A \otimes B := \{\, ab \mid a \in A,\; b \in B \,\}.
\end{equation}

The Minkowski sum was introduced in Minkowski's papers at the turn of the twentieth century, see, e.g., \cite{minkowski1897,minkowski1903}. 
The Minkowski product was studied more recently and requires a precise specification of the underlying multiplication. 
In this paper, we restrict ourselves to four-dimensional space $\mathbb{R}^4$, where the Minkowski product is naturally interpreted as the quaternionic product. 
Our aim is to apply Minkowski set operations to selected parametrically defined subsets of $\mathbb{R}^4$ and to visualize the results interactively using 4-D perspective and double orthogonal projections.

A modern definition of the Minkowski sum can be found, for example, in \cite[p.~126]{schneider1993}. 
The Minkowski sum has applications in various fields, e.g., motion planning \cite{lien2008} and geometric design \cite{muhlthaler2003}. 
Quaternionic Minkowski point set operations for unit quaternion sets are described and visualized using stereographic projection in \cite{farouki2019}. 
Related approaches based on the complex number product in $\mathbb{R}^2$ can be found, for example, in \cite{farouki2000,farouki2001,li2016}. 
Three-dimensional sets involving the wedge product are discussed in \cite{velichova2013a,velichova2013b}, and more generally in \cite{velichova2015,velichova2017}.

Quaternions have been closely connected with their geometric interpretation since their introduction by Hamilton \cite{hamilton1853,hamilton1866}. 
A comprehensive textbook on quaternions with an emphasis on computer visualization is \cite{hanson2006}. 
An overview of the use of the quaternionic product in $\mathbb{R}^4$ in kinematics is given in \cite{nawratil2016}. 
For recent applications to the construction of minimal surfaces from a polynomial perspective, see \cite{altavilla2026}.

Since the action of unit quaternions induces rotations and provides a pa\-ra\-me\-tri\-za\-tion of the 3-sphere, they are often visualized using conformal stereographic projection \cite{sanderson2018}. Minkowski quaternionic operations on circles on a unit 3-sphere in stereographic projection also proved to be essential for finding surfaces in $\mathbb{R}^3$ that contain two circles through each point \cite{skopenkov2019}.

In our approach, quaternions represent parametrically defined sets (curves, 2-surfaces, and 3-surfaces) in $\mathbb{R}^4$. 
For the demonstration of various properties, we therefore employ more intuitive methods for visualizing 4-D objects in 3-D: double orthogonal projection, a generalization of Monge's projection \cite{zamboj2018a,zamboj2018b}, and 4-D perspective \cite{rada2021,zamboj2025}. 
Our visualizations, created in \emph{Wolfram Mathematica}, are interactive 3-D models that often include dynamic manipulators for selected parameters.

\subsection{Minkowski quaternionic point set operations in $\mathbb{R}^4$}

We briefly summarize the necessary preliminaries on quaternions from the four-dimensional point of view adopted in this paper.

Let $A(a_0,a_1,a_2,a_3)$ and $B(b_0,b_1,b_2,b_3)$ be two quaternions considered as points in $\mathbb{R}^4$. 
The sum $A+B$ of the quaternions $A$ and $B$ is defined componentwise by

\begin{equation}
A+B=\bigl(a_0+b_0,\; a_1+b_1,\; a_2+b_2,\; a_3+b_3\bigr).
\end{equation}

The product $AB$ of the quaternions $A$ and $B$ is given by
\begin{equation}
\begin{split}
AB=\bigl(
a_0 b_0 - a_1 b_1 - a_2 b_2 - a_3 b_3,\;
a_0 b_1 + a_1 b_0 + a_2 b_3 - a_3 b_2,\;\\
a_0 b_2 - a_1 b_3 + a_2 b_0 + a_3 b_1,\;
a_0 b_3 + a_1 b_2 - a_2 b_1 + a_3 b_0
\bigr),
\end{split}
\end{equation}
or, equivalently, in scalar--vector form,
\begin{equation}
\begin{split}
AB
&=(a_0,\mathbf{a})(b_0,\mathbf{b}) \\
&=\left(
a_0 b_0 - \mathbf{a}\cdot\mathbf{b},\;
a_0 \mathbf{b} + b_0 \mathbf{a} + \mathbf{a}\times\mathbf{b}
\right),
\end{split}
\end{equation}
where $\mathbf{a}=(a_1,a_2,a_3)$ and $\mathbf{b}=(b_1,b_2,b_3)$ denote the vector parts of $A$ and $B$, respectively.

The neutral element of multiplication is the quaternion $1=(1,0,0,0)$.Since quaternion multiplication is non-commutative, we distinguish between 
left and right products. For quaternions $A$ and $B$, the left product of $B$ by $A$ is $AB$ and the right product of $B$ by $A$ is $BA$. Each quaternion $A$ has a conjugate 
\[
A^\ast=(a_0,-a_1,-a_2,-a_3).
\]
The norm of a quaternion $A$ is defined by
\[
|A|=\sqrt{AA^\ast}=\sqrt{a_0^2+a_1^2+a_2^2+a_3^2}.
\]
A quaternion $U$ satisfying $|U|=1$ is called a unit quaternion. For each quaternion $A$ with $|A|\neq 0$, the multiplicative inverse is given by
\[
A^{-1}=\frac{A^\ast}{|A|^2},
\] 
so that $A A^{-1}=1$. Because multiplication is non-commutative, for quaternions $A$ with $|A|\neq 0$ and $B$ we distinguish between the left division $A^{-1}B$ and the right division $BA^{-1}$ of $B$ by $A$.

Similarly to complex numbers, a quaternion $A$ admits a trigonometric form
\[
A = |A|\left(\cos\frac{\varphi}{2} 
+ \frac{\mathbf{a}}{|\mathbf{a}|}\,\sin\frac{\varphi}{2}\right),
\qquad \varphi \in [-\pi,\pi],
\]
provided that $\mathbf{a}\neq \mathbf{0}$.

From now on, we interpret a quaternion $A$ as the position vector of a point $A$ 
(keeping the same notation) in $\mathbb{R}^4$. The sum $A+B$ of two points coincides with the classical definition from analytic geometry in $\mathbb{R}^2$ or $\mathbb{R}^3$. 
The quaternionic product $AB$ can be interpreted as a rotation composed with a scaling of $B$ by $A$, or vice-versa, in $\mathbb{R}^4$, 
analogously to the complex product in the Gaussian plane $\mathbb{C} \cong \mathbb{R}^2$. 

Note that points $A$ and $B$ may also be regarded as 0-dimensional point sets in \eqref{eq:sum} and \eqref{eq:product}. 
Hence, $A+B = A\oplus B$ and $AB = A\otimes B$.

By introducing parameters, we extend our point sets from individual points to curves, surfaces, and higher-dimensional sets. 
In general, let $A(u_1,\dots,u_s)$ and $B(v_1,\dots,v_r)$ be parametrically defined point sets of dimensions $s$ and $r$, respectively, including the special case of $0$-dimensional point sets $A$ and $B$ (without parameters). 
Their Minkowski sum $A \oplus B$ and product $A \otimes B$ then form point sets of dimension $r+s$, provided that the embedding space has sufficiently high dimension. 

In this paper, we restrict ourselves to point sets of dimension at most three in $\mathbb{R}^4$. For example, consider two one-dimensional subsets $a(u)$ and $b(v)$, represented as curve segments with distinct parameters $u$ and $v$, given by their normalized parametric representations
\begin{equation}
\begin{split}
a(u) &= \bigl(a_0(u), a_1(u), a_2(u), a_3(u)\bigr), 
\qquad u \in [0,1],\\
b(v) &= \bigl(b_0(v), b_1(v), b_2(v), b_3(v)\bigr), 
\qquad v \in [0,1].
\end{split}
\end{equation}

The Minkowski sum $a(u)\oplus b(v)$ is a two-dimensional set — a surface patch in $\mathbb{R}^4$ — defined parametrically by $q(u,v)$ over the unit square $\lbrack0,1\rbrack \times \lbrack0,1\rbrack$:
\begin{equation}
q(u,v)=\bigl(
a_0(u)+b_0(v),\;
a_1(u)+b_1(v),\;
a_2(u)+b_2(v),\;
a_3(u)+b_3(v)
\bigr).
\end{equation}

The Minkowski difference $a(u)\ominus b(v)$ is defined analogously.

The Minkowski product $a(u)\otimes b(v)$ is again a two-dimensional set $p(u,v)$ defined over $\lbrack0,1\rbrack\times\lbrack0,1\rbrack$:
\begin{equation}
\begin{split}
p(u,v)=\bigl(
& a_0(u)b_0(v)-a_1(u)b_1(v)-a_2(u)b_2(v)-a_3(u)b_3(v),\\
& a_0(u)b_1(v)+a_1(u)b_0(v)+a_2(u)b_3(v)-a_3(u)b_2(v),\\
& a_0(u)b_2(v)-a_1(u)b_3(v)+a_2(u)b_0(v)+a_3(u)b_1(v),\\
& a_0(u)b_3(v)+a_1(u)b_2(v)-a_2(u)b_1(v)+a_3(u)b_0(v)
\bigr).
\end{split}
\end{equation}

Similarly, the Minkowski left and right divisions are given by 
$a(u)^{-1}\otimes b(v)$ and $b(v)\otimes a(u)^{-1}$, respectively.

If the point sets $A$ and $B$ share one or more parameters, these may be used to eliminate certain variables. 
In the case of multiplication (and so division), such a situation leads to a problem of factorization, which lies beyond the scope of this paper.

\subsection{4-D to 3-D methods of projection}

\begin{figure}
    \centering
    \includegraphics[width=.7\linewidth, trim= 0 50 0 50, clip]{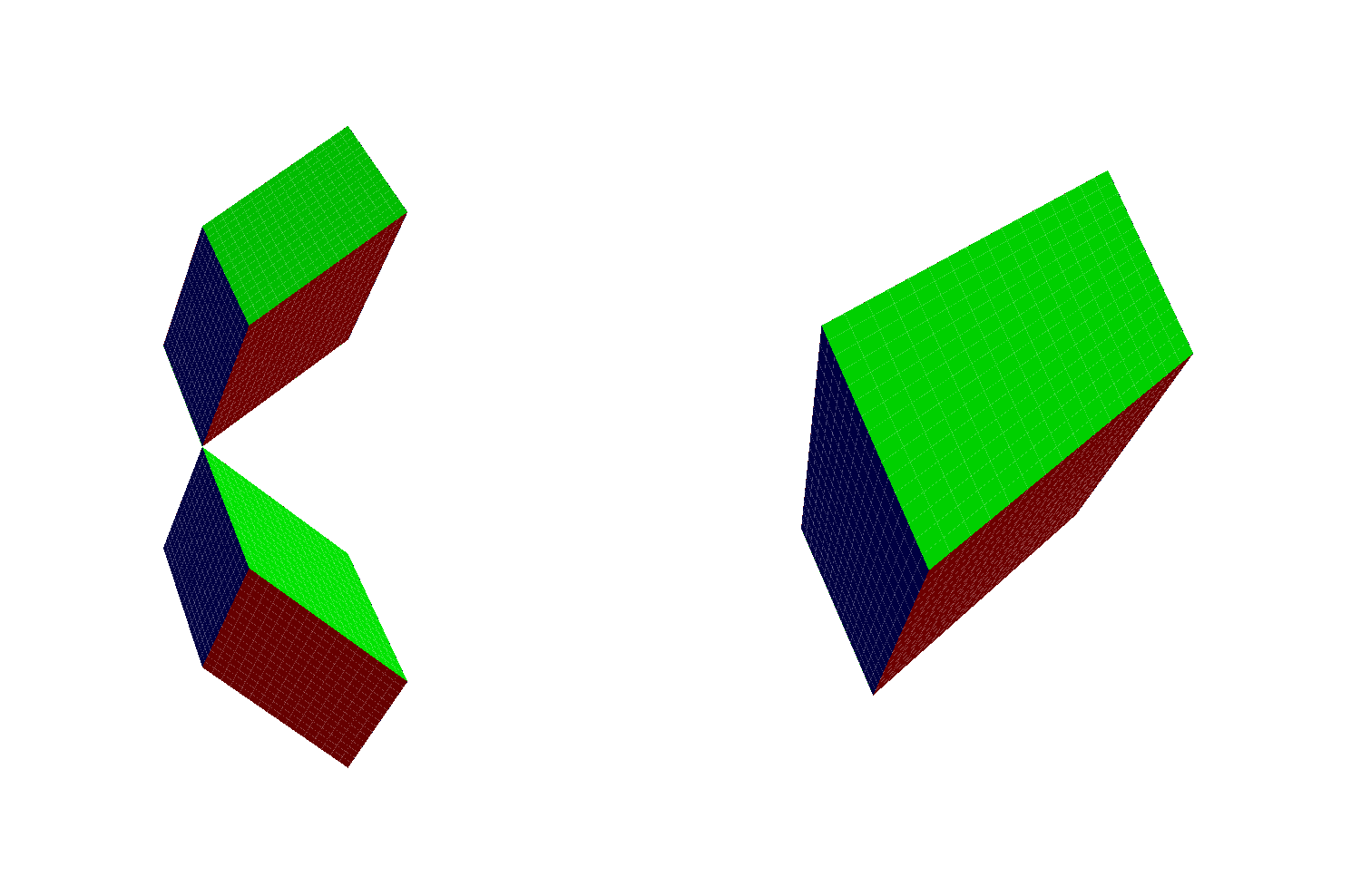}
    \caption{A 3-cube in DOP (left) and 4-D perspective (right).}
    \label{fig:cube}
\end{figure}

The visualizations presented throughout the paper are constructed using two projection methods from $\mathbb{R}^4$ to $\mathbb{R}^3$: double orthogonal projection of the 4-space $(x,y,z,w)$ onto two mutually perpendicular 3-spaces $(x,y,z)$ and $(x,y,w)$ (DOP), and central projection from a projection center in $\mathbb{R}^4$ onto a modeling 3-space $(x,y,w)$ (4-D perspective) (Figure~\ref{fig:cube}). The former is analogous to Monge's projection, while the latter corresponds to linear perspective generalized from 4-D to 3-D.

In the case of DOP, let $(x,y,w)$ denote the modeling 3-space, i.e., the space into which the point sets are projected and where the graphical visualizations are realized. 
Each point is orthogonally projected onto both 3-spaces $(x,y,z)$ and $(x,y,w)$. 
After projection, the 3-space $(x,y,z)$ is rotated about the common plane $(x,y)$ of the two 3-spaces into $(x,y,w)$ in such a way that the positive $z$-orientation becomes opposite to the positive $w$-orientation, so that their images overlap. The analytic representations of the projected coordinates of a point 
$A(a_0,a_1,a_2,a_3)$ are
\[
A_z(a_0,a_1,a_2)\quad \text{in } (x,y,z),
\qquad
A_w(a_0,a_1,a_3)\quad \text{in } (x,y,w).
\]

In the case of 4-D perspective, the modeling 3-space is again $(x,y,w)$. 
Each point in $\mathbb{R}^4$ is centrally projected from a projection center located at $(0,0,d,0)$ onto the modeling 3-space, where $d$ denotes the focal distance. 
The 4-D perspective maps a preimage point $A\in \mathbb{R}^4$ to its perspective image $A' \in \mathbb{R}^3$ by
\[
A(a_0,a_1,a_2,a_3) 
\longrightarrow 
A'\left(d\frac{a_0}{a_2},\, d\frac{a_1}{a_2},\, d\frac{a_3}{a_2}\right),
\qquad \text{for } a_2 \neq 0.
\]
In the case $a_2 = 0$, the image is given by
\[
A'(a_0,a_1,a_3).
\]

Since both projections reduce the dimension from four to three, certain geometric properties are inevitably lost. In this sense, the two methods complement each other. While the 4-D perspective provides a more global view and often a clearer intuitive impression within a single visualization, the double orthogonal projection preserves parallelism and, in special positions, also angles and distances.

\section{Selected examples}
The barycenter of this contribution is the geometric modeling of selected point sets in $\mathbb{R}^4$ using the Minkowski quaternionic sum and product.

\subsection{Clifford torus}

\begin{figure}[!htb]
    \centering
    \includegraphics[width=0.45\linewidth]{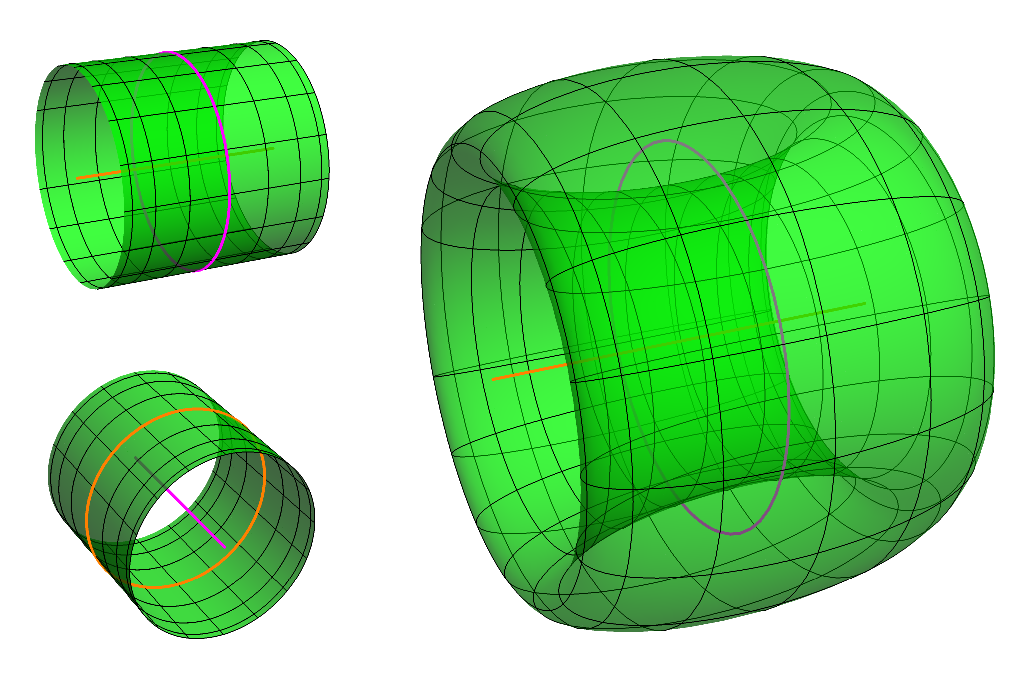}\hfill\includegraphics[width=0.45\linewidth]{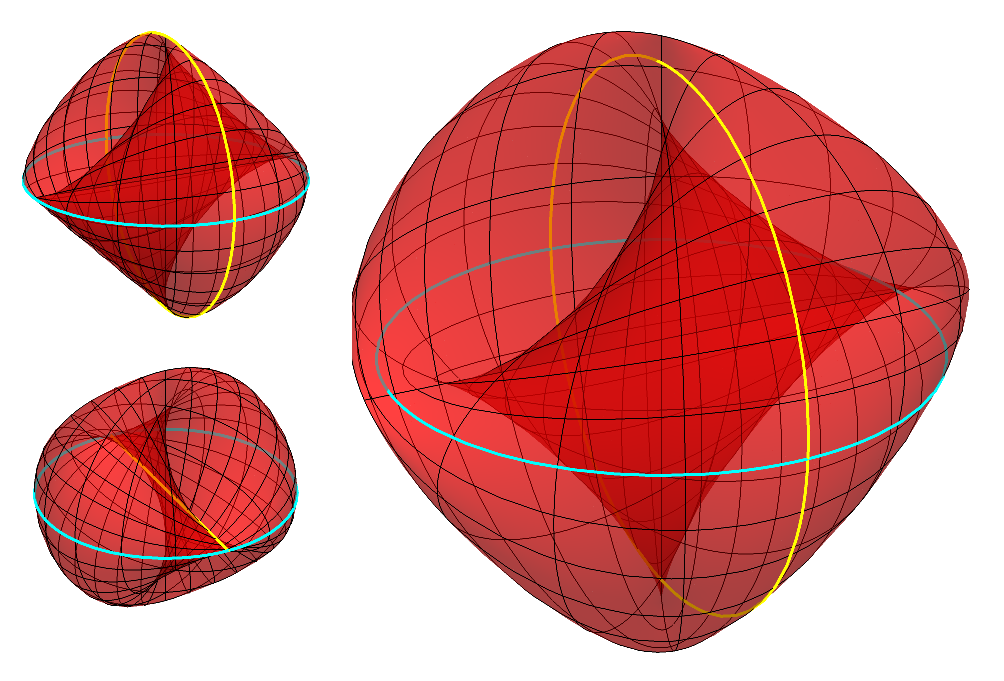}\\
    \includegraphics[width=0.15\linewidth]{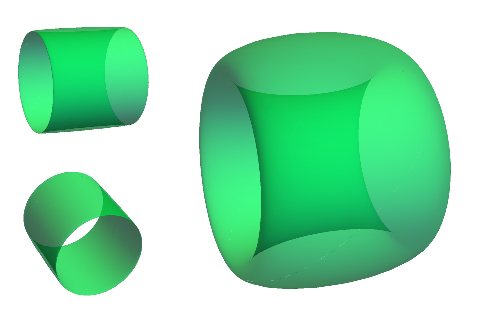} \includegraphics[width=0.15\linewidth]{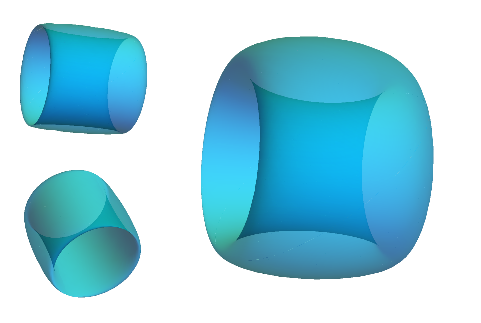} \includegraphics[width=0.15\linewidth]{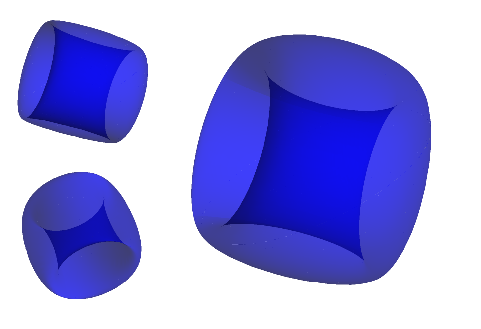} \includegraphics[width=0.15\linewidth]{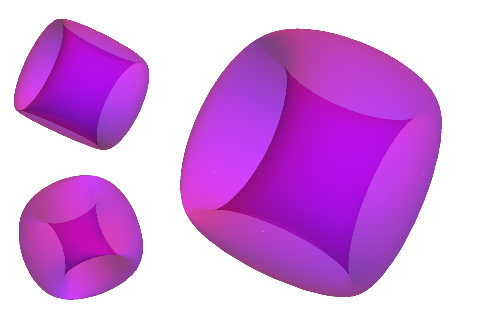} \includegraphics[width=0.15\linewidth]{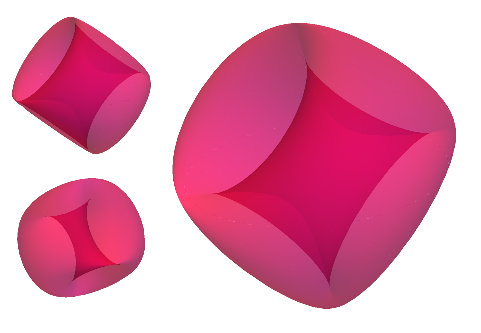}
    \caption{Clifford torus in DOP and 4-D perspective. The top left figure shows the Minkowski sum $c(u,v)$ generated according to Equations~\ref{eq:c-csum}, with the orange generating circle $c_1(u)$ and the purple generating circle $c_2(v)$. The top right figure shows the product $d(u,v)$ obtained from Equations~\ref{eq:c-cprod}, with generating circles $d_1(u)$ (cyan) and $d_2(v)$ (yellow). The bottom row depicts a sequence of rotations between $c(u,v)$ and $d(u,v)$.}

    \label{fig:c-c}
\end{figure}

We begin with the Minkowski sum of two circles:
\begin{equation}
\begin{split}
\label{eq:c-csum}
c_1(u) &= \left(-\frac{\sqrt{2}}{2}\cos u,\, 0,\, -\frac{\sqrt{2}}{2}\sin u,\, 0\right), 
\qquad u \in [0,2\pi],\\
c_2(v) &= \left(0,\, \frac{\sqrt{2}}{2}\cos v,\, 0,\, -\frac{\sqrt{2}}{2}\sin v\right), 
\qquad v \in [0,2\pi],\\
c(u,v) = c_1 \oplus c_2 
&= \left(
-\frac{\sqrt{2}}{2}\cos u,\,
\frac{\sqrt{2}}{2}\cos v,\,
-\frac{\sqrt{2}}{2}\sin u,\,
-\frac{\sqrt{2}}{2}\sin v
\right).
\end{split}
\end{equation}
The two-dimensional surface $c(u,v)$ is the well-known Clifford torus (Figure~\ref{fig:c-c}) lying on the 3-sphere (see \cite{rastanawi2022}), which will be examined in the following sections. 

For the Minkowski quaternionic product, we now consider a slightly different pair of circles:
\begin{equation}
\begin{split}
\label{eq:c-cprod}
d_1(u) &= (\cos u,\, \sin u,\, 0,\, 0), 
\qquad u \in [0,2\pi],\\
d_2(v) &= (0,\, \cos v,\, 0,\, \sin v), 
\qquad v \in [0,2\pi],\\
d(u,v) = d_1 \otimes d_2 
&= \bigl(
-\cos v \sin u,\,
\cos u \cos v,\,
-\sin u \sin v,\,
\cos u \sin v
\bigr).
\end{split}
\end{equation}

The two-dimensional surface $d$ is again a Clifford torus, positioned differently from $c$. 
The exact rotation mapping $d$ onto $c$ is obtained via the Minkowski quaternionic product
\begin{equation}
\begin{split}
    c(u_1,u_2)= 
\left(\frac{1}{2}, \frac{1}{2}, \frac{1}{2}, \frac{1}{2}\right)
\otimes d(v_1,v_2),\\
v_1=\frac{1}{2}(u_1+u_2),\;
v_2=\frac{1}{2}(u_1-u_2)+\frac{\pi}{4}.
\end{split}
\end{equation}
Observe that $d$ is multiplied from the right by a unit quaternion. 

\subsection{Quadratic cone}

\begin{figure}[!htb]
    \centering
    \includegraphics[width=\linewidth,trim=0 40 0 0, clip]{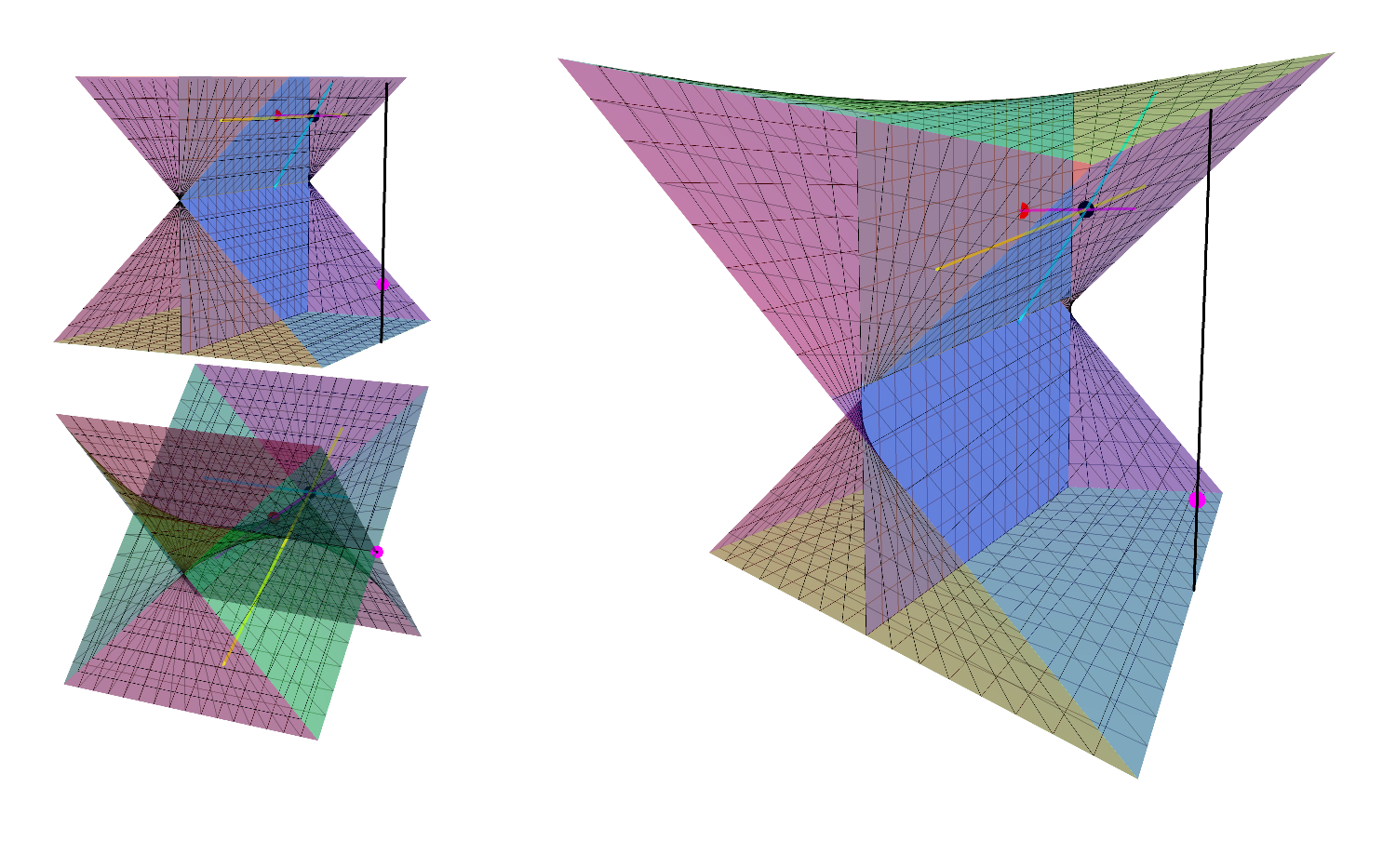}
   \caption{A portion of the quadratic cone $c$ given by Equations~\ref{eq:l-pprod}, parametrized over the cube $(u,v_1,v_2)\in[-1,1]\times[-1,1]\times[-1,1]$, visualized in DOP (left) and 4-D perspective (right). The black line represents $c_1(u)$, and the blue plane represents $c_2(v_1,v_2)$. The points on $c_1(u)$ and $c_2(v_1,v_2)$ correspond to the parameter values $u=-0.6$, $v_1=0.4$, and $v_2=0.5$. The black point is their product, and it lies on three rulings of $c(u,v_1,v_2)$.}

    \label{fig:l-p}
\end{figure}

For the next example, let the two point sets be a line and a plane. We construct their Minkowski quaternionic product.

\begin{equation}
\begin{split}
\label{eq:l-pprod}
c_1(u) &= (1,\,0,\,0,\,u), 
\qquad u \in \mathbb{R},\\
c_2(v_1,v_2) &= (0,\,v_1,\,0,\,v_2), 
\qquad (v_1,v_2) \in \mathbb{R} \times \mathbb{R},\\
c(u,v_1,v_2) = c_1 \otimes c_2 
&= (-u v_2,\, v_1,\, u v_1,\, v_2).
\end{split}
\end{equation}

Observe that $c(u,v_1,v_2)$ satisfies the simple implicit equation
\begin{equation}
xy + zw = 0,
\end{equation}
which defines a singular ruled hyperquadric in $\mathbb{R}^4$ with a singular point at $(0,0,0,0)$ and signature $(2,2)$. 

If we restrict the parameters to $u \in [-1,1]$ and $(v_1,v_2)\in [-1,1]\times[-1,1]$, 
we obtain a bounded portion of $c(u,v_1,v_2)$ (Figure~\ref{fig:l-p}). Its boundary consists of parts of hyperbolic paraboloids (for fixed $v_1$ or $v_2$) and planes (for fixed $u$).

Note that the intersection of $c$ with the unit 3-sphere is precisely the above-mentioned Clifford torus.

\begin{figure}[!htb]
    \centering
    \includegraphics[width=0.5\linewidth]{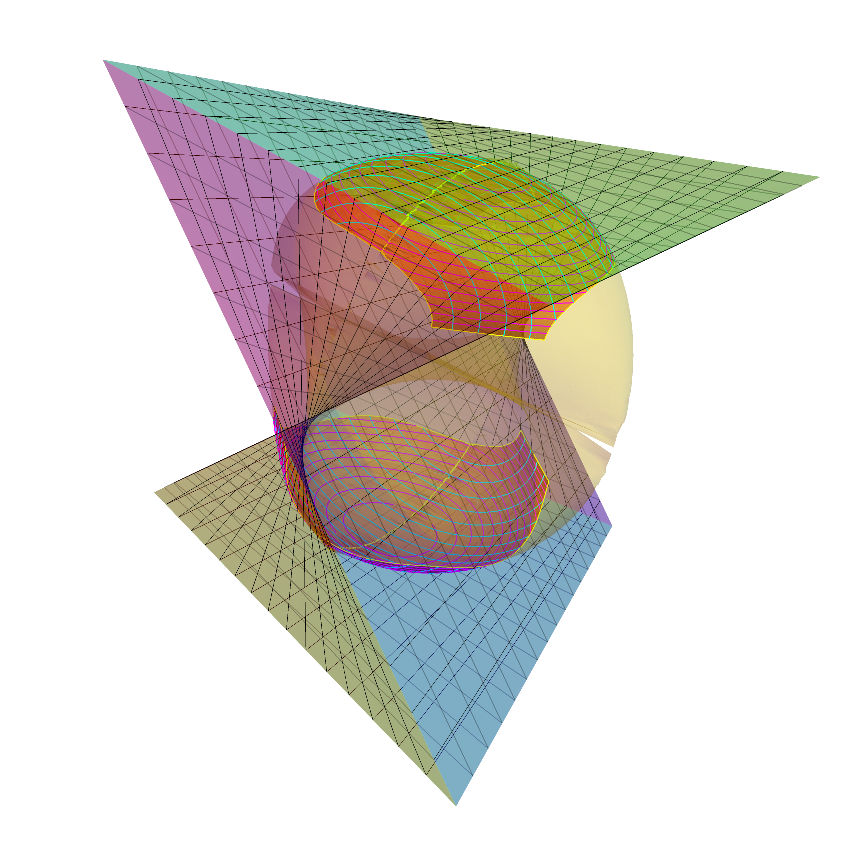}\\
    \includegraphics[width=0.3\linewidth]{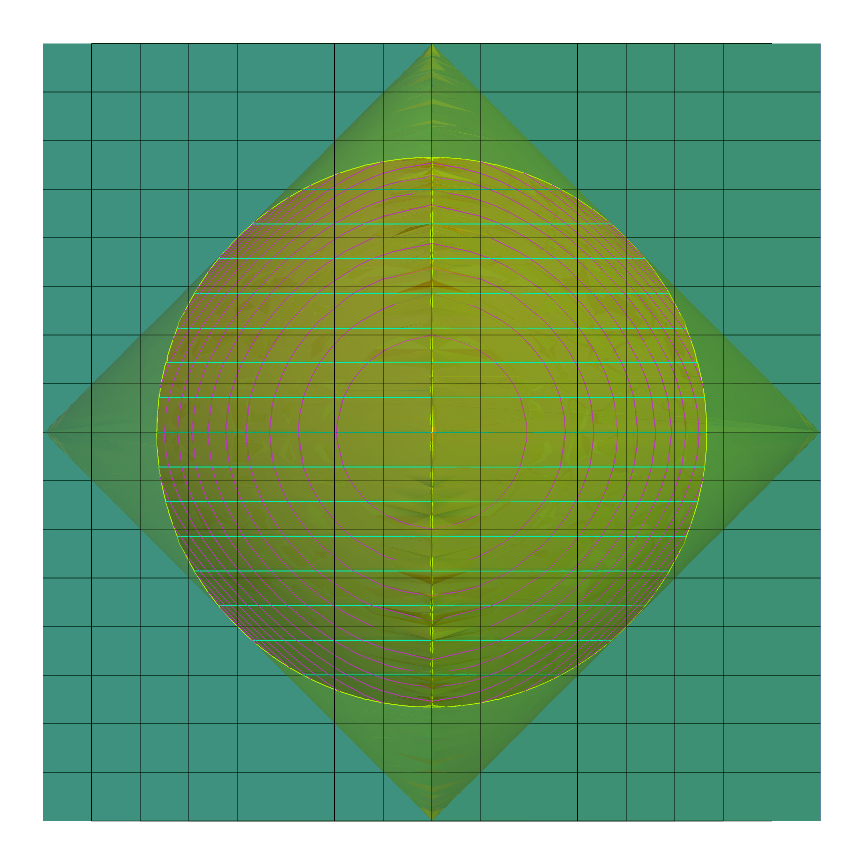}
    \includegraphics[width=0.3\linewidth]{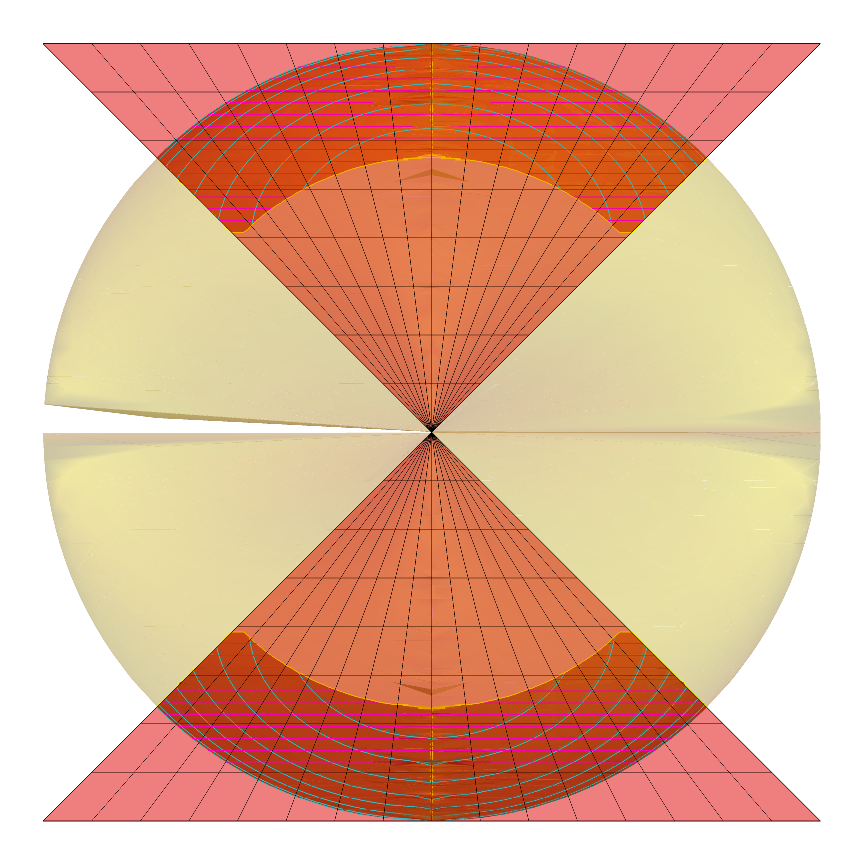}
    \includegraphics[width=0.3\linewidth]{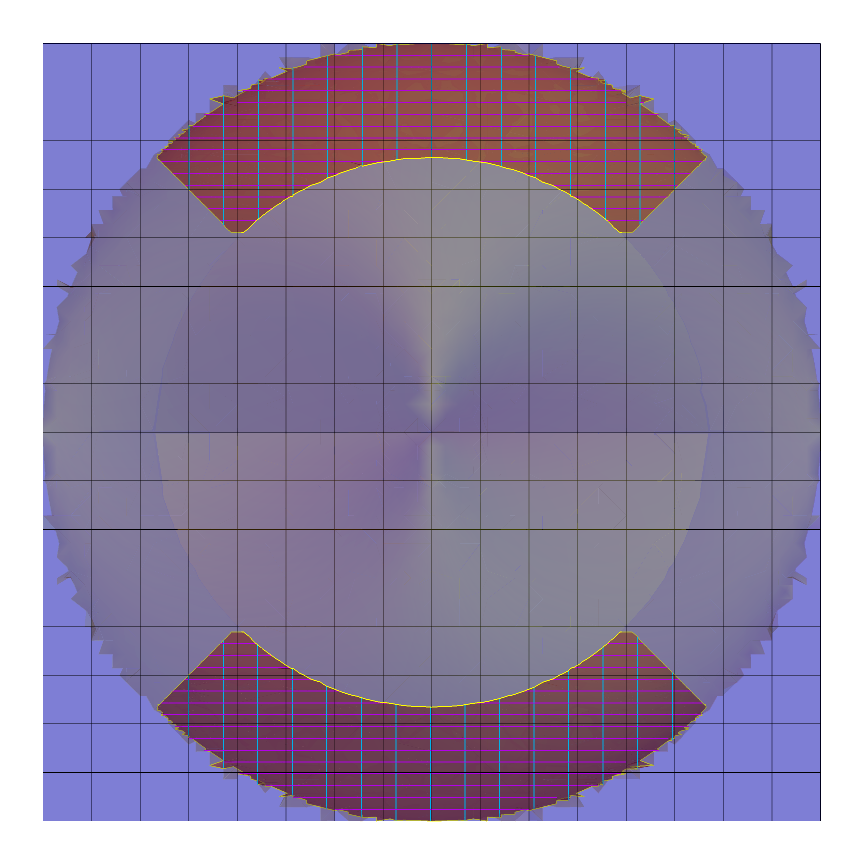}
    \caption{(Top) The intersection of the quadratic cone with the unit 3-sphere forms a Clifford torus. The highlighted region shows the portion of the torus corresponding to the restricted part of the cone. (Bottom) Orthogonal projections onto $(x,y)$, $(x,w)$, and $(y,w)$, respectively.}

    \label{fig:qcone-intersection}
    \end{figure}

\subsection{3-Sphere}

We have already referred several times to the unit 3-sphere. We now construct it as a Minkowski quaternionic product. Let the point sets be a circle and a 2-sphere:
\begin{equation}
\begin{split}
\label{eq:hopf}
c_1(u) &= (\cos u,\, \sin u,\, 0,\, 0), 
\qquad u \in [0,2\pi],\\
c_2(v_1,v_2) &= (\cos v_1,\, 0,\, \sin v_1 \cos v_2,\, \sin v_1 \sin v_2), 
\\ & \qquad (v_1,v_2) \in [0,\pi]\times[0,2\pi],\\
c(u,v_1,v_2) = c_1 \otimes c_2 
&=\\ =\bigl(
\cos u \cos v_1,\, &
\sin u \cos v_1,\,
\sin v_1 \cos (u+v_2),\,
\sin v_1 \sin (u+v_2)
\bigr).
\end{split}
\end{equation}

\begin{figure}[!htb]
    \centering
    \includegraphics[width=\linewidth, trim= 0 80 0 50, clip]{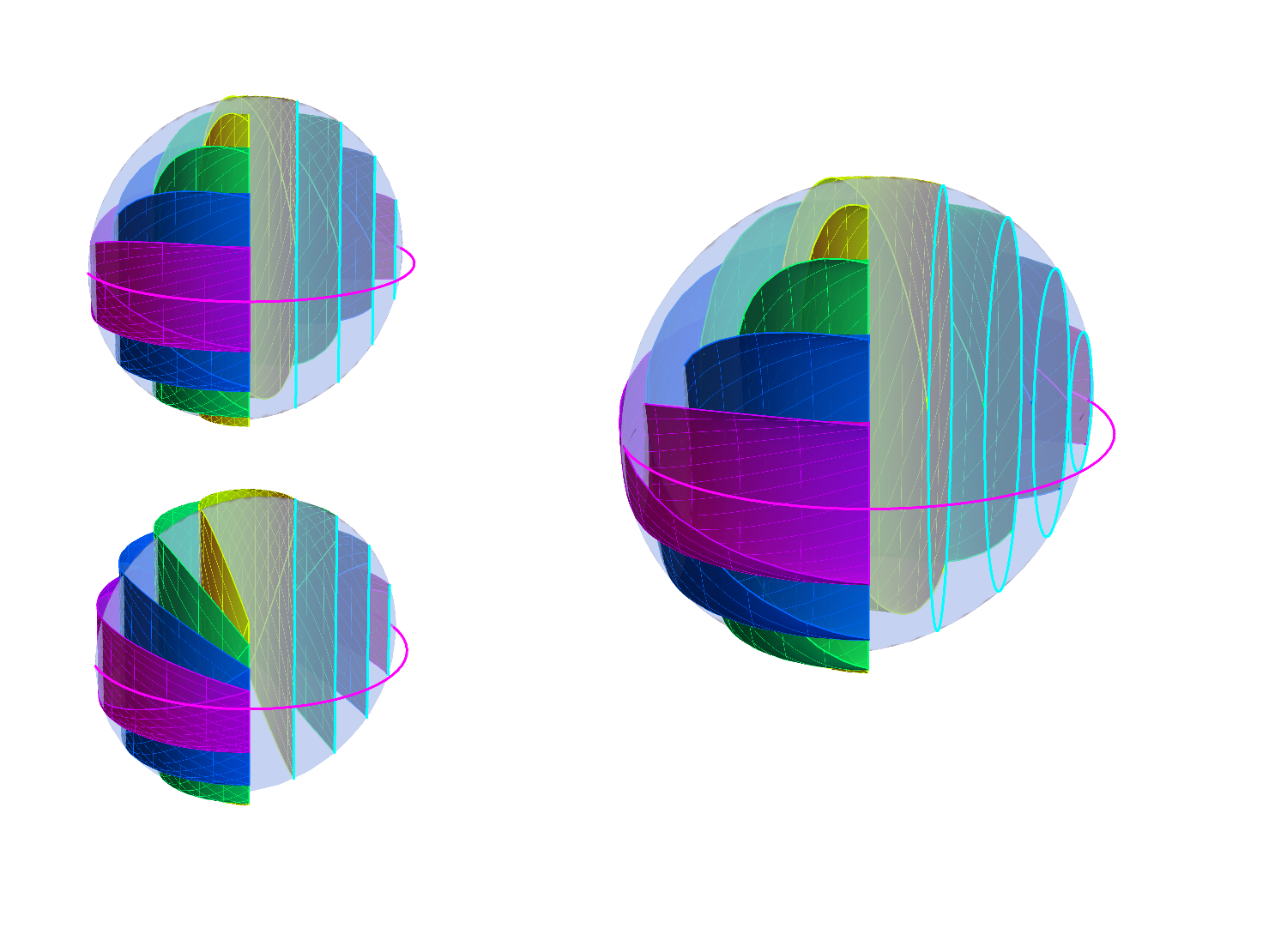}
    \caption{3-sphere generated by Equations~\ref{eq:hopf} as the product of the magenta circle $c_1(u)$ and light blue 2-sphere $c_2(v_1,v_2)$ (degenerated to a disk). The portions of tori refer to the cyan circles on the generating 2-sphere for selected choices of the parameter $v_1$.}
    \label{fig:hopf}
\end{figure}

In this case, the product $c(u,v_1,v_2)$ parametrizes the unit 3-sphere in Hopf coordinates with parameters $u$, $v_1$, and $v_2$. Figure~\ref{fig:hopf} illustrates a sequence of tori filling the 3-sphere. For more details on circles and tori on a 3-sphere see \cite{zamboj2021}. 

\subsection{Some more rulings}
We proceed with another example of geometric modeling using the Minkowski quaternionic product. Let the point sets be a line and a circle; subsequently, we will replace the circle by a helix.
\begin{equation}
\begin{split}
c_1(u) &= (1,\,-u,\,0,\,0), 
\qquad u \in \mathbb{R},\\
c_2(v) &= (0,\, \cos v,\, 0,\, \sin v), 
\qquad v \in [0,2\pi],\\
c(u,v) = c_1 \otimes c_2 
&= (u \cos v,\, \cos v,\, u \sin v,\, \sin v).
\end{split}
\end{equation}

\begin{figure}[!htb]
    \centering
    \includegraphics[width=.7\linewidth]{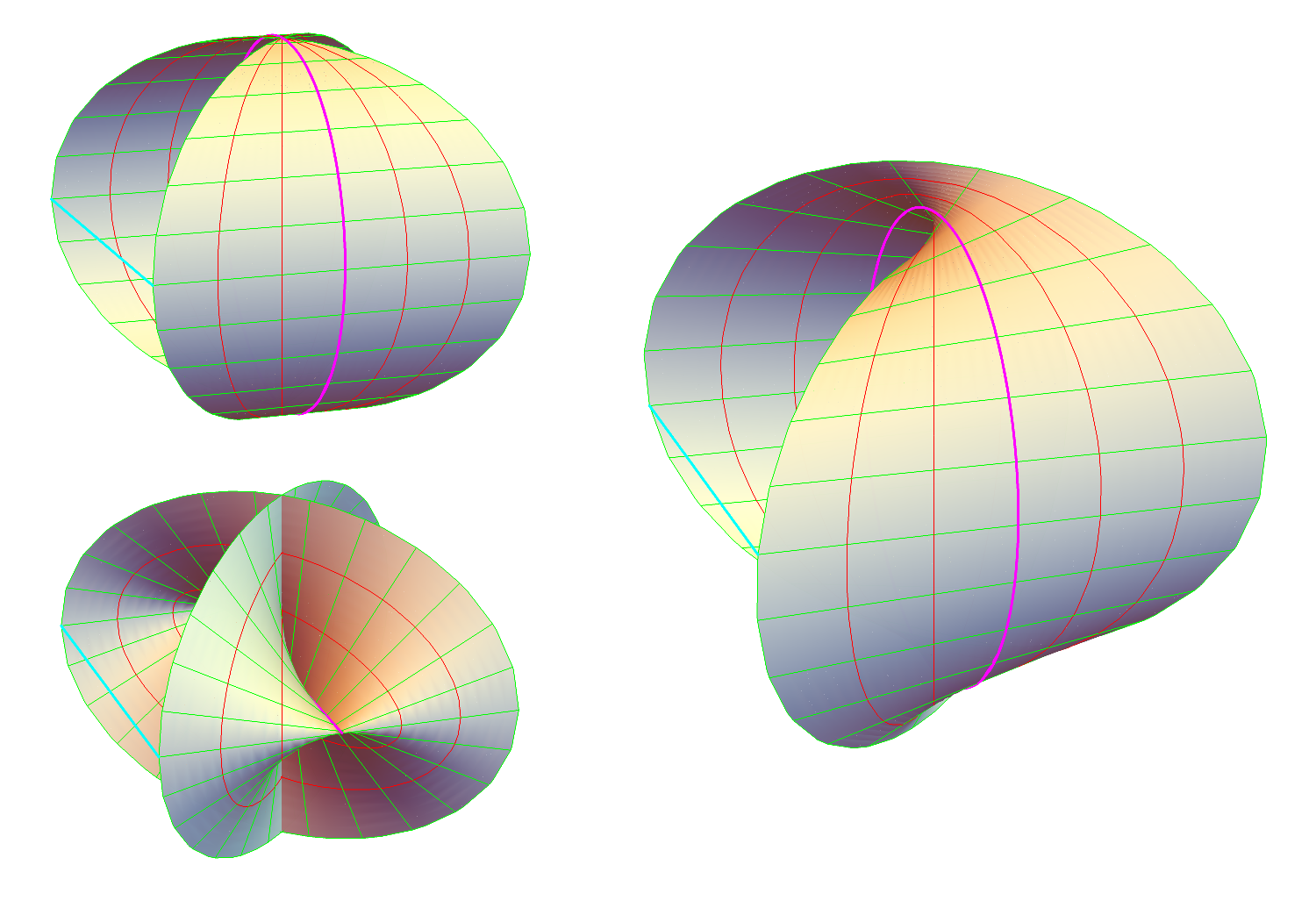}
    \caption{The 4-D perspective of the Minkowski product of the line segment $c_1(u)$ (cyan) and the circle $c_2(v)$ (magenta), consisting of families of lines and circles. The input parameters are $u\in[-1,1]$ and $v\in[-\pi,\pi]$. The orthogonal projection onto $w=0$ is Plücker's conoid.}

    \label{fig:line-circle}
\end{figure}

First observe (Figure~\ref{fig:line-circle}) that the orthogonal projection of $c(u,v)$ onto the hyperplane $w=0$ is the well-known surface called Plücker's conoid,
\[
(u \cos v,\, \cos v,\, u \sin v).
\]
For more details about this surface, see \cite{stachel2022}. 

Recall the kinematic interpretation of quaternion multiplication as a composition of a 4-D rotation and a homothety. Such transformations preserve lines and circles (although line segments and circle radii are scaled). Therefore, the product $c_1(u) \otimes c_2(v)$ can be interpreted in two equivalent ways: either the points on the line $c_1$ rotate and scale the circle (left product), or the points on the circle $c_2$ rotate and scale the line $c_1$ (right product). Consequently, the resulting two-dimensional surface $c(u,v) \subset \mathbb{R}^4$ contains both a family of circles and a family of lines.

\begin{figure}[!htb]
    \centering
    \includegraphics[width=.7\linewidth, trim=0 50 0 0,clip]{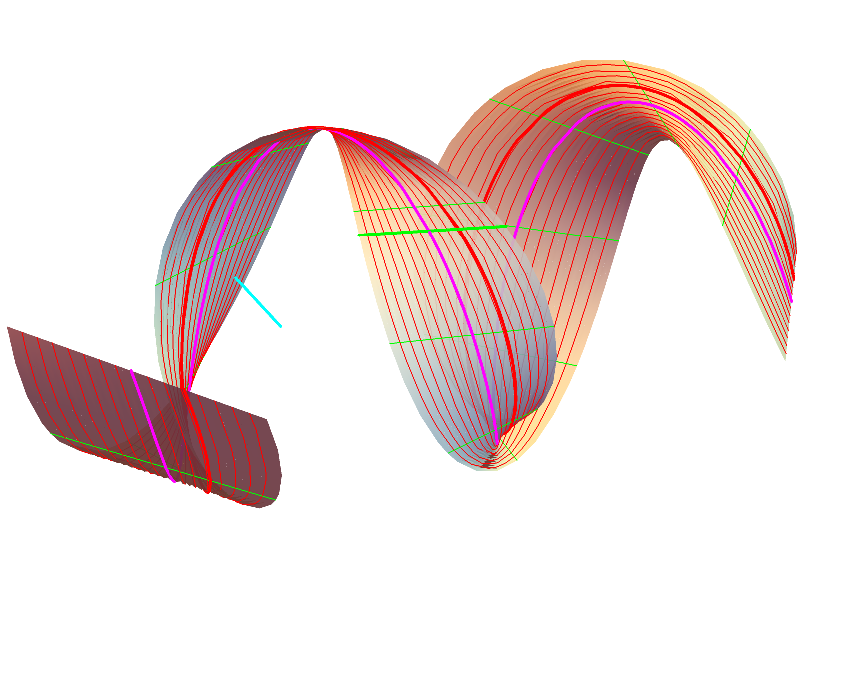}
    \caption{The 4-D perspective of the Minkowski product of the line segment $c_1(u)$ (cyan) and the helix $d_2(v)$ (magenta), consisting of families of lines and helices. The input parameters are $t=2$, $u\in[-0.5,0.5]$, and $v\in[-2\pi,2\pi]$.}

    \label{fig:line-helix}
\end{figure}

We now generalize the previous example by replacing the circle with a helix:
\begin{equation}
\begin{split}
\label{eq:line-helix}
c_1(u) &= (1,\,-u,\,0,\,0), 
\qquad u \in \mathbb{R},\\
d_2(v) &= \left(t\frac{v}{2\pi},\, \cos v,\, 0,\, \sin v\right), 
\qquad v \in \mathbb{R}, \quad t \in \mathbb{R} \text{ constant},\\
d(u,v) = c_1 \otimes d_2 
&= \left(
t\frac{v}{2\pi} + u\cos v,\,
- t u \frac{v}{2\pi} + \cos v ,\,
u \sin v,\,
\sin v
\right).
\end{split}
\end{equation}

Since a helix is not a closed curve (affinely in $\mathbb{R}^4$), the initial and final rulings do not coincide. The resulting ruled two-dimensional surface $d(u,v)$ also contains a family of helices; see Figure~\ref{fig:line-helix}.

\section{Conclusion and future work}

Minkowski point set operations establish a neat connection between algebraic and geometric representations. In particular, the choice of the product operation brings arithmetic laws into an intuitive geometric interplay. In this paper, we have adopted a geometric modeling perspective to construct and generalize selected objects in four-dimensional space. Finally, we have employed double orthogonal and perspective projections to highlight some of their characteristic properties.

\begin{figure}[!h]
    \centering
    \includegraphics[width=0.25\linewidth, trim= 0 20 0 20, clip]{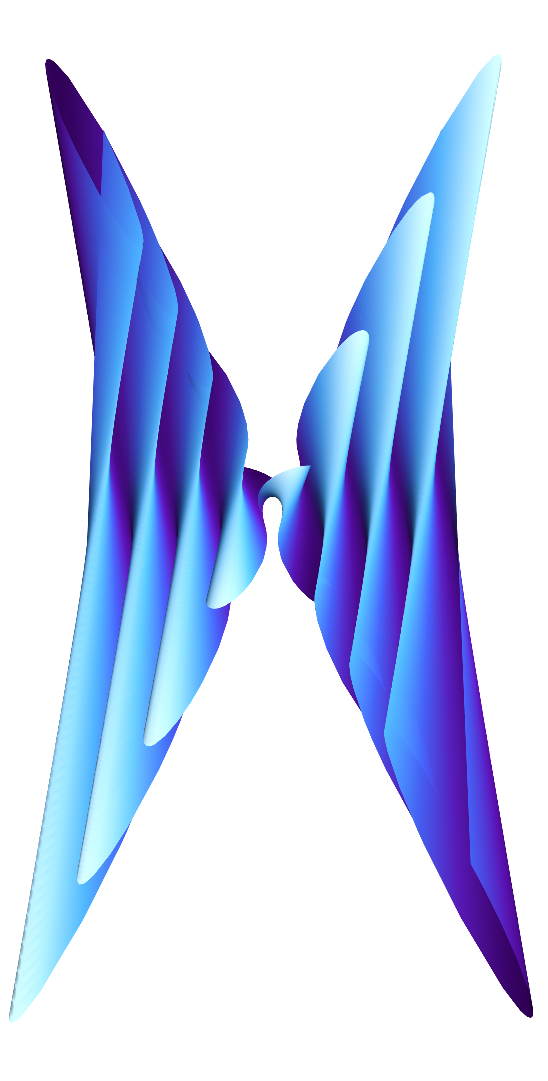}
  \caption{A butterfly generated by Equations~\ref{eq:line-helix}. The input parameters are $t=2\pi$, $u\in[-1,1]$, and $v\in[-8\pi,8\pi]$.}
  \label{fig:butterfly}
\end{figure}

Spectating the beauty, one can easily slip onto the path of excitement leading to infinite galleries of models. After all, see Figure~\ref{fig:butterfly}. 
There are many possible directions to proceed when supplementing the intuitive visual examination of Minkowski operations on point sets. 
Be it a classification of the point sets based on the generating sets (e.g., points, lines, circles, planes, spheres), pursuing properties such as sidedness, kinematic applications, or algebraic properties (factorization, roots, etc.).

%
%

\end{document}